\documentclass{article}

\usepackage{amsmath}
\usepackage{tikz}
\usepackage{amsfonts}
\usepackage{amsthm}
\usepackage{amssymb}
\usepackage{graphicx}
\usepackage{color}
\usepackage{verbatim} 
\usepackage{hyperref}
\usepackage{enumerate}
\usepackage{float}
\usepackage{tabulary}
\usepackage[utf8]{inputenc}
\theoremstyle{plain}

\theoremstyle{definition}

\newenvironment{packed_item}{
\begin{itemize}
  \setlength{\itemsep}{1pt}
  \setlength{\parskip}{0pt}
  \setlength{\parsep}{0pt}
}{\end{itemize}}

\title{Murder at the Asylum}
\author{Tanya Khovanova}
\date{}

\begin{document}
\maketitle

\begin{abstract}
I describe a puzzle I wrote for the 2018 MIT Mystery Hunt which introduced new types of people in logic puzzles. I discuss the puzzle itself, the solution, and the mathematics behind it.
\end{abstract}

\textit{In memory of Raymond Smullyan.}

I was a member of the team \textit{Death and Mayhem} which won the 2017 MIT Mystery Hunt. My team's reward was to write the 2018 MIT Mystery Hunt. By the way, some people consider it a punishment rather than a reward. As one colleague said to me on our win, ``Congratudolences!'' 

I decided to write a logic puzzle in tribute to Raymond Smullyan, who recently died. I ended up with a real monstrosity of a puzzle. 

In Section~\ref{sec:puzzle}, I present the puzzle with no explanatory comments. Section~\ref{sec:discussion} discusses the mathematical theory behind the puzzle. The solution is in Section~\ref{sec:solution}.

\section{Puzzle: Murder at the Asylum}\label{sec:puzzle}

Detective Radstein was scouting vacation spots and ended up on a peculiar island. Inhabitants of this island are very particular about the truth values of their statements. There are three types of people: \textit{truth-tellers}, \textit{liars} and \textit{alternators}. Truth-tellers only say things they believe are true; liars only say things they believe are false; and alternators switch at every statement between the two. The alternators switch according to their own internal counter. That means, if a stranger starts a conversation with the alternator, it is not known, \textit{a priori}, whether the alternator's first statement is true or false.

A murder was committed last night at the island's asylum, and Detective Radstein was invited to investigate.

Independently of being obsessive about their statements, there are three other types of people at the asylum: \textit{sane}, \textit{delusional} and \textit{partial}. Sane people believe true facts are true. The medical system is very efficient on the island: all the doctors are sane, and none of the patients is completely sane. 

The patients are either delusional or partial. A delusional patient only believes in facts that are not true. Moreover, they are very particular about their delusions. For example, if two plus two is four, the delusional patient believes that two plus two is not four. That is, they believe the exact negation of the true statement.

The partial patients combine traits of delusional patients and sane people. They switch the truthfulness of their beliefs after each statement they make. So a partial patient in the sane state can believe that two plus two is four, but after they make a statement, for example, that today is Wednesday, they switch to believing that two plus two is not equal to four. Similar to alternators, they switch their sanity according to an internal counter. That means, a stranger doesn't know, \textit{a priori}, whether the first statement is made in the sane or insane state.

When Radstein arrived at the asylum, he started by examining the victim's room. The victim was a patient named Zack who lived alone in his room. The victim was lying in his bed, his skull crushed. Next to the victim was a bloody bronze statue of a ballerina. It was obviously the murder weapon. Radstein put on his gloves and picked up the statue. It was very heavy: he could barely lift it. It would require significant strength to use this statue as a murder weapon. Radstein examined the victim. The man obviously died in his sleep. Radstein estimated the time of death as somewhere between 2am and 5am. 

Radstein looked around the room: all the windows were locked. The only way anyone could have entered the room was through the door, but the door lock had not been meddled with. Radstein also found a statue stand in the common room.

The asylum is small. There are only nine people there: Ann, Beth, Cedric, David, Eve, Fiona, Grace, Holly, and Ian. This is not counting the dead guy. It is worth noting that the members of the asylum do not talk about things they do not know about.

The detective went to a room where all the asylum folks gathered. They immediately started talking.

\begin{packed_item}
\item Ann: Zack was my lover.
\item Beth: Zack was my lover.
\item Cedric: Zack was my lover.
\item David: Zack was my lover.
\item Eve: Zack was my lover.
\item Fiona: Zack was my lover.
\item Grace: Zack was my lover.
\item Holly: Zack was my lover.
\item Ian: Zack was my lover.
\end{packed_item}

``That was amusing,'' thought Radstein. He decided to control the conversation: ``Yes or no, are you a patient?''

\begin{packed_item}
\item Ann: Yes.
\item Beth: Yes.
\item Cedric: No.
\item David: Yes.
\item Eve: Yes.
\item Fiona: Yes.
\item Grace: Yes.
\item Holly: No.
\item Ian: Yes.
\end{packed_item}

``With all these strange alternators and partials, why don't I ask the same question again,'' thought Radstein. Here are the answers.

\begin{packed_item}
\item Ann: Yes.
\item Beth: Yes.
\item Cedric: Yes.
\item David: No.
\item Eve: No.
\item Fiona: Yes.
\item Grace: Yes.
\item Holly: Yes.
\item Ian: No.
\end{packed_item}

``Hmm, will this help me find alternators and partials?'' thought Radstein, ``What happens if an alternator is also a partial? There is also the curious issue of their delusions. Is there a way to distinguish between insane truth-tellers and sane liars? What do insane people believe that they believe? Let me try something else.'' He asked: Yes or no, do you believe that you are a patient?''

\begin{packed_item}
\item Ann: Yes.
\item Beth: No.
\item Cedric: No.
\item David: No.
\item Eve: Yes.
\item Fiona: No.
\item Grace: No.
\item Holly: No.
\item Ian: No.
\end{packed_item}

While Radstein tried to collect his thoughts, the suspects got excited and started to discuss their love lives again.

\begin{packed_item}
\item Ann: Beth was Zack's lover.
\item Beth: I believe Cedric was Zack's lover.
\item Cedric: David was Zack's lover.
\item David: I believe Eve wasn't Zack's lover.
\item Eve: Fiona was Zack's lover.
\item Fiona: I believe that Grace wasn't Zack's lover.
\item Grace: I believe Holly wasn't Zack's lover.
\item Holly: Ian was Zack's lover.
\item Ian: I believe Ann wasn't Zack's lover.
\end{packed_item}

Then they moved on to the murder.

\begin{packed_item}
\item Ann: I believe the patients' doors are locked from 11pm till 7am and only doctors have keys to the patients' rooms.
\item Beth: I believe there is a killer among the doctors.
\item Cedric: I believe I saw a person carrying the statue to the victim's room at 2am.
\item David: I believe that the person who unlocked the door is not a doctor or not a lover.
\item Eve: I believe there exists a man who is not innocent.
\item Fiona: I believe there is a killer among the delusional patients.
\item Grace: I believe there is a non-lover who is not innocent.
\item Holly: I believe that only Beth, Cedric, David, and Ian are strong enough to lift the statue.
\item Ian: I believe that the person who brought the statue to Zack's room is not a lover.
\end{packed_item}

Radstein came back to his hotel room to write a report. At the end of the report he wants to have a list for each person. In addition to solving the crime Radstein wants to figure out who is who. It is not clear if the person who unlocked the room and/or the person who brought the statue into the room also committed the murder. It is curious that for guilt too he needs three categories just in case.

The first bullet on the list is related to sanity. It will say whether the person is 
\begin{packed_item}
\item partial,
\item delusional, or
\item sane.
\end{packed_item}
The second bullet is about truthfulness. It will say whether the person is 
\begin{packed_item}
\item an alternator, 
\item a liar, or
\item a truth-teller.
\end{packed_item}
The third bullet is related to guilt: 
\begin{packed_item}
\item an accomplice,
\item guilty of murder, or
\item innocent.
\end{packed_item}

Please, help Detective Radstein.

\section{Discussion}\label{sec:discussion}

Raymond Smullyan wrote numerous logic books with many beautiful logic problems. He started with two kinds of people: \textit{knights} and \textit{knaves}. Knights always tell the truth, knaves always lie. Raymond Smullyan coined the names for knights and knaves in his 1978 book \textit{What Is the Name of This Book?} \cite{RSWhat}.  His problems often took places on an island to guarantee that the main characters of the problems were of particular type.

Raymond Smullyan introduced sane and insane people in his book \textit{The Lady or the Tiger?} in 1982 \cite{RSLOrT}. Sane people believe true facts are true, while insane people believe that false facts are true. A chapter of Smullyan's book is devoted to different asylums, in which each member is either a doctor or a patient. All the doctors are supposed to be sane and all the patients insane. In addition, they all are truthful. Inspector Craig is sent to inspect the asylums. Here is one of the problems in the book:

\begin{quote}
What statement can a sane patient make to prove that they should be released?
\end{quote}

The answer: ``I am a patient.'' Indeed, sane doctors believe they are doctors and wouldn't say that they are patients. Similarly, insane patients believe that they are not patients and can't say they are patients.

One might think that insane truth-tellers are equivalent to sane liars. Indeed, when talking about facts they both lie. Here is another problem from that book that shows how interesting the concept of insanity is:

\begin{quote}
Inspector Craig asked a person in the asylum, ``Are you a patient?'' The reply was, ``I believe so.'' Is there something wrong with the asylum?
\end{quote}

This is a very tricky question. If a sane person says ``I believe A'', then s/he believes A. What happens with an insane person? If an insane person says s/he believes A, it means the negation of the statement must be true. It means, they really believe NOT A. For example, an insane person believes that two plus two is not four. But being insane s/he believes that s/he believes that two plus two IS four. In other words, an insane truth-teller will say, ``I believe that two plus two is four,'' which is a lie as it should be. That means, when we look at statements about beliefs, the underlying statement is true for all truth-tellers. It follows that, in the problem, the person at the asylum is a patient, but we do not have a way to know whether this is a sane or an insane patient.

What would happen if insane people can also lie? Here is an example that shows a pitfall that we avoided in the MIT Mystery Hunt puzzle. If a liar believes that two plus two is five, then s/he can lie that two plus two is six. S/he can also lie that two plus two is four. As a result, we get either a true or false statement, which doesn't provide us with any conclusion. This is why this puzzle was so specific about beliefs. An insane person in this puzzle believes that two plus two is not four. Hence, an insane liar always says that two plus two is four.

We already showed that truth-tellers, both sane and insane, would say, ``I believe two plus two is four.'' What about liars? Sane liars believe that two plus two is four. They also believe that they believe that two plus two is four. Being liars they would say that they believe that two plus two is not four. Similarly, insane liars believe that two plus two is not four. They also believe that they believe that two plus two IS four. Being liars they would say that they believe that two plus two is not four. When talking about beliefs, both sane and insane liars behave the same way. They would both say. ``I believe not A,'' where A is a true fact.

Let us combine what we learned into Table~\ref{tab:SITL}. The first column describes types of people, the second column is the answer to the question, ``Are you a patient?'' The third column is the answer to the question, ``Do you believe that you are a patient?'' We use Y for ``yes,'' and N for ``no.'' Remember that doctors are sane and patients are insane.

\begin{table}[htb]\label{tab:SITL}
\begin{center}
\begin{tabular}{| r | r | r |}
\hline
sane truth-teller & N & N\\
sane liar & Y & N \\
insane truth-teller & N & Y \\
insane liar & Y & Y \\
\hline
\end{tabular}
\end{center}
\caption{Sane/insane truth-tellers/liars answer questions.}
\end{table}

We see that these two questions allow us to differentiate between the types of people that do not switch their truthfulness or their sanity.

But this is not enough, as this puzzle has more types of people: people who can switch their truthfulness and their sanity. Alternators switch at every statement between telling the truth and lying. The patients are either delusional or partial. A delusional patient is always insane. Partial patients switch between sane or insane states.

The only person who thinks he is a patient is a partial in the sane state. Here is what we can deduce about the asylum patients and doctors from the following statements:

\begin{itemize}
\item ``I am a patient.'' This can be sane lying, delusional lying, insane partial lying, and sane partial telling the truth.
\item ``I am not a patient.'' This can be sane telling the truth, delusional telling the truth, sane partial lying, and insane partial telling the truth.
\end{itemize}

Patients in insane or sane states believe that they believe that they are patients. Doctors believe that they believe that they themselves are not patients. We can analyze the following statements:

\begin{itemize}
\item ``I believe I am a patient.'' This can be sane lying, delusional telling the truth, sane partial telling the truth, and insane partial telling the truth.
\item ``I believe I am not a patient.'' This can be sane telling the truth, delusional lying, sane partial lying, and insane partial lying.
\end{itemize}

How do we differentiate between sane people who might be truth-tellers, liars, or alternators? The big idea is to ask the same question twice: ``Is two plus two equal to four?'' The truth-teller will answer Y and Y, liar N and N, and the alternator will switch between Y and N depending on the starting state. Repeating this question, we not only can find the alternator but also the alternator together with their starting state.

Let us denote the sanity of people prior to the first question as follows: S is sane, D is delusional, Pi partial insane, and Ps partial sane. Similarly for truthfulness, let us denote T for truth-tellers, L for liars, At for alternators telling the truth, and Al for alternators lying.
The big idea to differentiate between all these types, is to ask four questions in a row:

\begin{itemize}
\item Are you a patient?
\item Are you a patient?
\item Do you believe that you are a patient?
\item Do you believe that you are a patient?
\end{itemize}

Table~\ref{tab:all} shows how each type answers the four questions.

\begin{table}[htb]\label{tab:all}
\begin{center}
\begin{tabular}{| r | r | r | r | r | r | r | r | r | r | r | r | r | r | r | r |}
\hline
ST & SL & SAt & SAl & DT & DL & DAt & DAl & PiT & PiL & PiAt & PiAl & PsT & PsL & PsAt & PsAl \\ 
\hline
N & Y & N & Y &  N & Y & N & Y &  N & Y & N & Y &  Y & N & Y & N\\
N & Y & Y & N &  N & Y & Y & N &  Y & N & N & Y &  N & Y & Y & N\\
N & Y & N & Y &  Y & N & Y & N &  Y & N & Y & N &  Y & N & Y & N\\
N & Y & Y & N &  Y & N & N & Y &  Y & N & N & Y &  Y & N & N & Y\\
\hline
\end{tabular}
\end{center}
\caption{Question answers per type.}
\end{table}

By asking these four questions in a row, we can differentiate between types. Now we are ready to solve the puzzle.

\section{Solution}\label{sec:solution}

Let's skip the lovers round in which every person claimed to be Zach's lover, We will count that round of statements as round zero, and will ignore it for now. 

Now let's look at the next three rounds of questions, which we will call rounds one, two, and three. We also mark the individuals' sanity/truthfulness states before round one.

We saw how the four questions can differentiate between the types. Being evil, I put only three out of the four questions in the puzzle. So we can't completely differentiate. Each sequence of answers produces two potential types. By analyzing three consecutive questions about patients we get the following possibilities summarized in Table~\ref{tab:types}.

\begin{table}[htb]\label{tab:types}
\begin{center}
\begin{tabular}{| r | r | r |}
\hline
NNN & ST, PsAl & \\
YYY & SL, PsAt & Ann \\
NYN & SAt, PsL & Cedric, Holly \\
YNY & SAl, PsT & Eve \\
YYN &  DL, PiAl & Beth, Fiona, Grace \\
NNY &  DT, PiAt &  \\ 
YNN &  DAl, PiL & David, Ian \\
NYY &  DAt, PiT & \\
\hline
\end{tabular}
\end{center}
\caption{Suspects' types.}
\end{table}

We do not know who is who yet. But we know something. Namely, we can deduce information from some of the statements. Let us look at factual statements (not beliefs). Both DL and PiAl always tell the truth about simple facts. That means that Beth, Fiona, and Grace are Zack's lovers. Now that we know some lovers we can start unraveling who is who from Round 4 (where they talk about each others' love lives), working backwards. We start with Ann, who talks about Beth, and we already know that Beth is a lover.

\begin{packed_item}
\item Ann: Beth was Zack's lover. Therefore, Ann is PsAt in insane lying state, rather than SL. Ann is a lover.
\item Ian: I believe Ann wasn't Zack's lover. Therefore, Ian is PiL in sane state, rather than DAl in true state. Ian is not a lover.
\item Holly: Ian was Zack's lover. Therefore, Holly is SAt in lying state, rather than PsL in insane state. Holly is not a lover.
\item Grace: I believe Holly wasn't Zack's lover. Therefore, Grace is PiAl in sane true state rather than DL. Grace is a lover.
\item Fiona: I believe that Grace wasn't Zack's lover. Therefore, Fiona is DL rather than PiAl in sane true state. Fiona is a lover.
\item Eve: Fiona was Zack's lover. Therefore, Eve is SAl in true state, rather than PsT in insane state. Eve is a lover.
\item David: I believe Eve wasn't Zack's lover. Therefore, David is PiL in sane state, rather than DAl in true state. David isn't a lover.
\item Cedric: David was Zack's lover. Therefore, Cedric is SAt in lying state, rather than PsL in insane state. Cedric isn't a lover.
\item Beth: I believe Cedric was Zack's lover. Therefore, Beth is DL, rather than PiAl in sane true state. Beth is a lover.
\end{packed_item}

In the last round the asylum folks only talked about their beliefs. That means we only need to know if they are telling the truth in order to decipher the true statements:

\begin{packed_item}
\item The patients' doors are locked from 11pm till 7am and only doctors have keys to the patients' rooms.
\item There are no killers among the doctors.
\item There was a person carrying the statue to the victim's room at 2am.
\item The person who unlocked the door is a doctor and a lover.
\item The men are innocent.
\item There are no killers among the delusional patients.
\item All non-lovers are innocent.
\item Only Beth, Cedric, David, and Ian are strong enough to lift the statue.
\item The person who brought the statue to Zack's room is a lover.
\end{packed_item}

The person who unlocked the door is a doctor and a lover. This could only be Eve. The person who carried the statue to the room must be either Beth, Cedric, David, or Ian. We also know that the men are innocent. This leaves Beth as the only person who could have carried the statue. We also know that the killer might only be a partial: this leaves Ann, David, Grace, or Ian. Given that the men are innocent, we are left with Ann and Grace. Given that none of them alone could have lifted the statue, both of them are killers.

This ends the logic part of the puzzle. But in puzzle hunts, the answer is a word or a phrase. The last step, where the hunters find this word, is called \textit{the extraction step}. So, what do we do with all this information?

Did you notice that in Radstein's report there are three categories: truthfulness, sanity, and guilt. There are also three types of people in each category. The frequent use of three, hints at the ternary extraction. Plus, the characteristics in the report do not match the order they were used in the puzzle. That indicates that this order is important. 

Did you notice that the names of all the suspects start with distinct letters from the beginning of the alphabet? This is a hint that we need to put the suspects in alphabetical order. Then for each suspect we assign a value based on Radstein's list. For every category, we use digits 0, 1, and 2 for the first, second, and third type on the list correspondingly. Starting with zero makes it easier to write the answers as ternary numbers. Then we convert ternary to decimals and decimals to letters. The process is summarized in Table~\ref{tab:results}.

\begin{table}[htb]\label{tab:results}
\begin{center}
\begin{tabular}{| r | r | r | c |}
\hline
Ann 		& PsAt, killer		& 001, 1  & A \\
Beth 		& DL, accomplice	& 110, 12 & L \\
Cedric 		& SAt, innocent 	& 202, 20 & T \\
David 		& PiL, innocent		& 012, 5  & E \\
Eve		& SAl, accomplice 	& 200, 18 & R \\
Fiona 		& DL, innocent 		& 112, 14 & N \\
Grace 		& PiAl, killer		& 001, 1  & A \\
Holly 		& SAt, innocent		& 202, 20 & T \\ 
Ian		& PiL, innocent		& 012, 5  & E \\
\hline
\end{tabular}
\end{center}
\caption{The extraction.}
\end{table}

The answer to the puzzle is \textbf{ALTERNATE}.


\begin{thebibliography}{9}


\bibitem{RSWhat} Raymond Smullyan, \textit{What is the Name of this Book?}, Prentice-Hall, (1978)

\bibitem{RSLOrT} Raymond Smullyan, \textit{The Lady or the Tiger?}, Knopf, (1982)

\end{thebibliography}
\end{document}